\documentclass[12pt]{article}

\usepackage{amsfonts}
\usepackage{amssymb}
\usepackage{enumerate}

\newtheorem{st}      {Theorem}
\newtheorem{prop}  {Proposition}
\newtheorem{lem} {Lemma}

\newtheorem{rmk} {Remark}[section]
\newtheorem{rmks} {Remarks}[section]

\newcommand{\f}{\ensuremath{\varphi}}
\newcommand{\rmap}{\longrightarrow}
\newcommand{\Boxe}{\raisebox{.8ex}{\framebox}}
\newcommand{\lmap}{\longleftarrow}

\newcommand{\U}{\ensuremath{\mathcal{U}}}

\newcommand{\Ha}{\ensuremath{H}}
\newcommand{\G}{\ensuremath{G}}

\newcommand{\nV}{\ensuremath{\mathcal{V}}}

\newcommand{\R}{\ensuremath{\mathcal{R}}}

\newcommand{\el}{\ensuremath{\mathcal{L}}}
\newcommand{\F}{\ensuremath{\mathcal{F}}}
\newcommand{\es}{\ensuremath{\mathcal{S}}}

\newcommand{\D}{\ensuremath{\mathcal{D}}}

\renewcommand{\Im}{\mathop{\mathrm{Im}}}
\newcommand{\Hol}{\mathop{\mathit{Hol}}\nolimits}
\newcommand{\hol}{\mathop{\mathrm{hol}}\nolimits}
\newcommand{\Id}{\mathop{\mathrm{Id}}\nolimits}
\newcommand{\Ind}{\mathop{\mathrm{Ind}}\nolimits}
\newcommand{\Mon}{\mathop{\mathit{Mon}}\nolimits}
\newcommand{\Top}{\mathop{\mathit{Top}}\nolimits}
\newcommand{\Trace}{\mathop{\mathrm{Tr}}\nolimits}
\newcommand{\Chern}{\mathop{\mathrm{Ch}}\nolimits}
\newcommand{\Etale}{\mathop{\mathit{Etale}}\nolimits}
\newcommand{\Hochs}{\mathop{\mathrm{HH}}\nolimits}
\newcommand{\rank}{\mathop{\mathrm{rank}}}
\newcommand{\supp}{\mathop{\mathrm{supp}}}

\newcommand{\basic}{\mathrm{bas}}

\newcommand{\ps}{{\raise 1pt\hbox{\tiny (}}}

\newcommand{\pss}{{\raise 1pt\hbox{\tiny [}}}
\newcommand{\pdd}{{\raise 1pt\hbox{\tiny ]}}}
\newcommand{\pd}{{\raise 1pt\hbox{\tiny )}}}

\newcommand{\nG}[1]{\ensuremath{\G_{#1}}}
\newcommand{\nK}[1]{\ensuremath{K_{#1}}}

\def\cross{\ltimes}

\newcommand{\nH}[1]{\ensuremath{H_{\, #1\, }}}

\def\compose{{\raise 1pt\hbox{$\scriptscriptstyle\circ$}}}
\def\dcross{{\raise 0.5pt\hbox{$\scriptscriptstyle\boxtime$}}}

\advance\textheight by \topskip
\usepackage[all,tips]{xy}
\usepackage[all]{xy}
\UseComputerModernTips

\begin{document}

\title{Foliation groupoids and their cyclic homology\thanks{Research 
supported by NWO}}
\author {Marius Crainic and Ieke Moerdijk}
\pagestyle{myheadings}
\date {Department of Mathematics, Utrecht University, The Netherlands}
\maketitle



\section*{Introduction}

%
The purpose of this paper is to prove two theorems which concern
the position of \'etale groupoids among general smooth (or ``Lie'')
groupoids. Our motivation comes from the non-commutative geometry
and algebraic topology concerning leaf spaces of foliations. Here,
one is concerned with invariants of the holonomy groupoid of a
foliation \cite{CoOp, Wi}, such as the cohomology of its
classifying space \cite{Haefl}, the cyclic homology of its smooth
convolution algebra \cite{BrNi, Cra}, or the $K$-theory of the
$C^*$-convolution algebras. Many results here depend on the fact
that such a holonomy groupoid can be ``reduced'' to what is called a
complete transversal of the foliation, giving rise to an
equivalent \'etale groupoid. For \'etale groupoids (sometimes called
$r$-discrete groupoids in the literature \cite{Pat, Re}), the cyclic
homology, sheaf theory and classifying spaces are each well
understood, as is the relation between these.

Our first theorem provides a criterion for determining 
whether a given Lie
groupoid is equivalent to an \'etale one. We prove that this is the 
case if and only if all
the isotropy groups of the groupoid are discrete, or equivalently, 
exactly when the anchor
map of the associate Lie algebroid is injective. These conditions are 
often easy
to check in examples. 

We recall that the Lie algebroid of a Lie groupoid is an
infinitesimal structure which plays the same role as the Lie
algebra of a Lie group. Lie algebroids with injective anchor map
are the same things as foliations, so another way of phrasing our
first theorem is by saying that a Lie groupoid is equivalent to an
\'etale one, exactly when it integrates a foliation. For this
reason, we have decided to refer to these  groupoids as ``foliation
groupoids''. It is not a surprise to see that much of the standard
literature on foliations deals with foliation groupoids; for instance,
an overall assumption in \cite{MoSo} is the discreteness of the isotropy
groups. Our first theorem can also be seen as a general
``slice theorem'', which generalizes the reduction to transversals
for foliations and the slice theorem for infinitesimally free
actions of compact Lie groups. This slice theorem is expected to be
a special case of a more general slice theorem  conjectured by A.
Weinstein. We also prove that, among the Lie groupoids which
integrate a given foliation, the holonomy and monodromy groups are
extreme examples. Results of this kind, but formulated in terms of 
microdifferentiable
groupoids, go back to \cite{Pra, Brown, Mon}.

Our second theorem concerns the invariance of cylic type
homologies under equivalence. We prove that equivalent foliation
groupoids have isomorphic Hochschild, cyclic and periodic cyclic
homology groups. This invariance is perhaps not really surprising,
especially since analogous results for \'etale groupoids \cite{CrMo, Mrc} ,
and for the K-theory of C*-algebras associated to groupoids \cite{MRW, HiSk}
are well known (see also \cite{Cun}). Nonetheless, we believe our 
second theorem has
some relevance. The theorem implies that the cyclic type homologies 
of leaf spaces
are totally independent of the particular model of the holonomy
groupoid, and
its proof provides explicit isomorphisms (summarized in the Remark at the end).
The theorem also completes the computation for algebras associated to 
Lie group actions
with discrete stabilizers. Moreover, this
second theorem may in fact be an intermediate step toward a
similar result for (more) general Lie groupoids. (Observe in this
context that some parts of the proof, such as the $H$-unitality of the
convolution algebra, apply to general Lie groupoids.)

The plan of this paper is as follows. In the first 
section we have collected the preliminary definitions concerning Lie
groupoids, their Lie algebroids, and their cyclic homology. In the 
second section
we state the main results. Since our motivation partly came from a 
better understanding of (the
relation between different approaches to) the longitudinal index theorem for
foliations (see \cite{CoMo, Co3, Ni3, HeLa}), we have added a few 
brief comments
at the end of this section. Section \ref{sec3} contains the proof of 
the first theorem
and the related results, and Section \ref{sec4} contains the proof of 
the second theorem.
We also mention that the part concerned with the theory of Lie groupoids
(namely Theorem \ref{charact} and Proposition \ref{classific}, and 
their proofs in
Section \ref{sec3}) can be read independently of the
preliminaries on cyclic homology in Section \ref{sec1}.

\section{Preliminaries}
\label{sec1}
%
%
We begin by recalling the necessary definitions and notation
concerning groupoids and cyclic homology. Standard references
include \cite{McK, Pat, Haefl} for groupoids, and \cite{Co3, Ka, Lo} for
cyclic homology.

\paragraph{Groupoids:}\label{groupoids}%
A {\it groupoid} $\G$ is a (small) category in which every arrow is
invertible. We will write $\nG{0}$ and $\nG{1}$ for the set of
objects and the set of arrows in $\G$, respectively.
The source and target maps are denoted by $s, t: \nG{1}\rmap \nG{0}$,
while $m(g, h)= g\compose h$ is the composition, and $i(g)= g^{-1}$ 
denotes the inverse of $g$.
One calls $\G$ a  {\it smooth groupoid} if $\nG{0}$ and
$\nG{1}$ are smooth manifolds,
all the structure maps are smooth, and $s$ and $t$ are submersions.
Basic examples include Lie groups, manifolds, crossed products of
manifolds by Lie
groups, the holonomy and the monodromy groupoids of a foliation,
Haefliger's groupoid $\Gamma^q$, and groupoids associated to
orbifolds.

If $\G$ is a smooth groupoid and $X, Y \subset \nG{0}$,
we write $\G_{X}= s^{-1}(X)$, $\G^{Y}= t^{-1}(Y)$, 
$\G_{X}^{Y}= s^{-1}(X) \cap t^{-1}(Y)$. Note that $G_{X}^{X}$ has the 
structure of a groupoid
(the restriction of $\G$ to $X$). When $X= Y= \{ x\}$, $x\in \nG{0}$,
we simplify the notations to $G_{x}$, $G^{x}$, $G_{x}^{x}$;
these are submanifolds of $\nG{1}$, and $G_{x}^{x}$ is a Lie group,
called the \emph{isotropy group} of $G$ at~$x$.

The tangent spaces at $1_x$ of
$G_{x}$ form a bundle $\mathfrak{g}$ over $\nG{0}$, of
``$s$-vertical'' tangent vectors on $\nG{1}$; it is the restriction
along $u: \nG{0} \hookrightarrow \nG{1}$ 
of the vector bundle $T^s\bigl(G_1\bigr) =
\ker\bigl(ds: T\nG{1} \rightarrow T\nG{0}\bigr)$. The differential
$dt: T\nG{1} \rightarrow T\nG{0}$ of the target map induces a map
of vector bundles over $\nG{0}$,
\[ \alpha : \mathfrak{g} \rightarrow T\nG{0},\]
called the \emph{anchor map}.  Moreover, the space of sections
$\Gamma\mathfrak{g}$ is equipped with a Lie bracket
$[ \cdot \, , \cdot ]$. This bracket makes $\alpha$ into a Lie algebra
homomorphism $\Gamma(\alpha):\Gamma \mathfrak{g} \rightarrow
\mathcal{X}\bigl(\nG{0}\bigr)$ into the vectorfields on $\nG{0}$, satisfying
the identity $[X, fY] = f[X,Y]+ \alpha(X)(f) \cdot Y$ for any
$X, Y \in \Gamma \mathfrak{g}$ and
$f \in C^{\infty}\bigl(\nG{0}\bigr)$.  This structure
\[ \bigl(\mathfrak{g}, \, [\cdot \, , \cdot ] \, , \alpha\bigr) \]
is called {\it the Lie algebroid} of $G$, and briefly denoted
$\mathfrak{g}$ in this paper.

A {\it
homomorphism} $\varphi : \G \rightarrow \Ha$ between two smooth
groupoids is a smooth functor. Thus, it is given by two smooth
maps (both) denoted $\varphi : \nG{0} \rightarrow \nH{0}$ and
$\varphi: \nG{1} \rightarrow \nH{1}$, commuting with all the
structure maps ($\varphi\circ s = s \circ \varphi$, etc.). Such a
homomorphism is called an {\it essential equivalence} if the map $s
\pi_{2}: \nK{0} \times_{\nG{0}} \nG{1} \rmap \nG{0}$, defined on
the space of pairs $\ps y, g\pd \in \nK{0} \times \nG{1}$ with
$t\ps g\pd= \f \ps y\pd$, is a surjective submersion, and the
square
\begin{equation}\label{eq1} \xymatrix {
\nK{1} \ar[r]^-{\f_{1}} \ar[d]_-{\ps s, t\pd} & \nG{1}
\ar[d]^-{\ps s,t \pd}\\
\nK{0}\times \nK{0} \ar[r]^-{\f_{0} \times \f_{0}} & \nG{0}\times
\nG{0} }
\end{equation}
is a pullback. Two groupoids $\G_i$ are said to be {\it Morita
equivalent} if there exists a third groupoid $\G$, and essential
equivalences $\f_{i}:\G\rmap \G_i$ as above ($i\in \{1, 2\}$).

If $f: X\rmap \nG{0}$ is a smooth map, one defines
the pullback of $\G$ along $f$ as the groupoid $f^*(\G)$ whose
space of objects is $X$, and whose arrows between $x, y\in X$ are
the arrows of $\G$ between $f(x)$ and $f(y)$. When the map
$s\pi_2: X \times_{\nG{0}} \nG{1} \rmap \nG{0}$ is a surjective
submersion, the groupoid $f^*(\G)$ is smooth and the obvious
smooth functor $f^*(\G)\rmap \G$ is a Morita equivalence. For
instance, given a family $\U= \{ U_i\}$ of opens in
$\nG{0}$, we define the groupoid $\G_{\U}$ as the pullback along
$f:\coprod_i U_i \rmap \nG{0}$. If $\U$ is a covering, then
$\G_{\U}$ is Morita equivalent to $\G$. Also,
if $\G= \Hol(M, \F)$ is the holonomy groupoid of a
foliation $(M, \F)$, and $i_{T}: T\rmap M$ is a
{\it transversal} for $\F$ (recall that this means that $T$ intersects each
leaf transversally), 
then $i_{T}^{*}\bigl(\Hol(M, \F)\bigr)=\Hol_{T}(M, \F)$ is the reduced
holonomy groupoid of $\F$. If $T$ is a {\it complete transversal} 
(i.e. intersects
each leaf at least once), then $\Hol_{T}(M, \F)$ is the standard 
\'etale groupoid
(see below) which is Morita equivalent to $\Hol(M, \F)$.\\

A smooth groupoid $\G$ is called  {\it \'etale} (or 
$r$-discrete) if
the source map $s: \nG{1} \rmap \nG{0}$ is a local diffeomorphism.
This implies that all other structure maps are also local
diffeomorphisms. Basic examples are discrete groups, manifolds,
crossed products of manifolds by (discrete) groups, the reduced
holonomy groupoid of a foliation, Haefliger's groupoid $\Gamma^q$,
groupoids associated to orbifolds.

The category
$\Etale$ of \'etale groupoids (with generalized homomorphisms)
plays an essential role in the study of leaf spaces of foliations.
It should be viewed as an enlargement of the category of smooth
manifolds
\begin{equation}\label{inclcat} \Top\ \subset\ \Etale \
\end{equation}
to which many of the classical constructions from algebraic
topology extend: homotopy, sheaves, cohomology, compactly
supported cohomology, Leray spectral sequences, Poincar\'e
duality, principal bundles, characteristic classes etc. See 
\cite{CrMo, CraMoe, Haefl, Moe, Mrc}.

In extending these constructions, one often uses the
following property, typical of \'etale groupoids. Any arrow $g:
x\rmap y$ induces a (canonical) germ $\sigma_{g}: (U,
x)\tilde{\rmap} (V, y)$ from a neighborhood $U$ of $x$ in $\nG{0}$
to a neighborhood $V$ of $y$. Indeed, we can define $\sigma_{g} =
t \compose \sigma$, where $x\in U \subset \nG{0}$ is so small that
$s: \nG{1}\rmap \nG{0}$ has a section $\sigma: U\rmap \nG{1}$ with
$\sigma\ps x \pd = g$.\\

\paragraph{Convolution algebras and cyclic homology:} 
Let $\G$ be a
smooth groupoid. To define its smooth convolution algebra
$C_{c}^{\infty}\bigl(\G\bigr)$, one uses the convolution product, defined
for functions $\phi$, $\psi$ on $\G$ and $g\in\nG{1}$, by
\begin{equation}\label{convprod}
(\phi*\psi)(g)= \int_{g_1g_2= g} \phi(g_1)\psi(g_2)\;.
\end{equation}

We assume for simplicity that $\G$ is Hausdorff. (For 
general groupoids,
possibly non-Hausdorff, the construction of the convolution algebra
is slightly more involved \cite{CrMoHa}.) If $\G$ is \'etale, then 
the integration is simply
summation, but, in general, one has to give a precise meaning to
the integration in the previous formula. For this, some choices
have to be made. If one wants to work with complex-valued
functions $\phi, \psi \in C_{c}^{\infty}\bigl(\nG{1}\bigr)$, then one has to
fix a smooth Haar system for $\G$ (we refer to \cite{Re} for
precise definitions). Instead, it is possible to use a line bundle
$\el$ of ``densities'' which is isomorphic to the trivial bundle
(in a non-canonical way), and to work with compactly supported
smooth sections of $\el$, $\phi, \psi \in C_{c}^{\infty}\bigl(\G;
\el\bigr)$. Fixing a trivialisation of $\el$ induces a Haar system on
$\G$, and gives an isomorphism $C_{c}^{\infty}\bigl(\G; \el\bigr) \cong
C_{c}^{\infty}\bigl(\G\bigr)$.

Let us
recall Connes' choice of $\el$ \cite{CoOp}. Let $\mathfrak{g}$ be the Lie
algebroid of $\G$. Denote by $\D^{1/2}$ the line bundle on
$\nG{0}$ consisting of transversal half-densities. Writing $p= 
\dim(\mathfrak{g})$, the fiber of
$\D^{1/2}$  over $x\in \nG{0}$ consists of maps $\rho$ from the exterior power
$\Lambda^{p}\mathfrak{g}_{x}$ to $
\mathbb{C}$ such that $\rho(\lambda v)= |\lambda|^{1/2}\rho(v)$
for all $\lambda\in  \mathbb{R}$,
$v\in\Lambda^{p}\mathfrak{g}_{x}$. There is a similar bundle  $\D^{r}$ for any
$r$. The bundle of densities ($r=1$) is usually denoted by $\D$.
We put $\el= t^*\D^{1/2}\otimes s^{*}\D^{1/2}$. Then
(\ref{convprod}) makes sense for  $\phi, \psi \in
C_{c}^{\infty}\bigl(\G; \el\bigr)$. Indeed, looking at the variable $g_2=
h$, one has to integrate $\phi(gh^{-1})\psi(h)$ $\in
\D_{x}^{1/2}\otimes \D_{z}\otimes \D_{y}^{1/2}$ with respect to $z
\stackrel{h}{\lmap} y$ varying in $\G_{y}$. But $\D_{z}$ is
canonically isomorphic to the fiber at $h$ of the bundle of
densities on the manifold $\G_{y}$, hence the integration makes
sense and gives an element $(\phi*\psi)(g)\in \D_{x}^{1/2}\otimes
\D_{y}^{1/2}= \el_{g}$. 
In the sequel we will
omit $\el$ from the notation $C_{c}^{\infty}\bigl(\G; \el\bigr)$.

Given an algebra $A$, recall the definition of
Connes' cyclic complex $C_{*}^{\lambda}\bigl(A\bigr)$, and of Hochschild's
complex $C_{*}\bigl(A\bigr)$. 
The latter has $C_{n}\bigl(A\bigr)= A^{\otimes (n+1)}$, 
with boundary $b$ given by
\begin{eqnarray}
b(a_0, a_1, \, . \, .\, .\, , a_n)     & = &  b' (a_0, a_1, \, . \, 
.\, .\, , a_n) +
(-1)^n(a_{n}a_0, a_1, \, . \, .\, .\, , a_{n-1})\, ,\nonumber\\
b'(a_0, a_1, \, . \, .\, .\, , a_n) & = &  \sum_{i=0}^{n-1} (-1)^i(a_0,
\, . \, .\, .\, , a_ia_{i+1}, \, . \, .\, .\, , a_n) \, ,       \nonumber
\end{eqnarray}
while the cyclic complex is the quotient 
$C_{n}^{\lambda}\bigl(A\bigr):= A^{\otimes (n + 1)}/\Im(1- \tau)$ 
with boundary induced by $b$. 
Here $\tau$ is the signed cyclic
permutation:
\[ \tau(a_0, a_1, \, . \, .\, .\, , a_n)= (-1)^n (a_n, a_0, \, . \, 
.\, .\, ,  a_{n-1})\ .\]

Recall that the cyclic homology groups $HC_{*}\bigl(A\bigr)$ of $A$ are
computed by the complex $C_{*}^{\lambda}\bigl(A\bigr)$. Also, the Hochschild
homology groups $\Hochs_{*}\bigl(A\bigr)$ are computed by $C_{*}\bigl(A\bigr)$, 
provided
$A$ is {\it $H$-unital}. Recall that $H$-unitality  means that
$\bigl(C_*(A), b'\bigr)$ is acyclic, and it plays a crucial role in the
excison theorems for cyclic homology \cite{Wod}.
For instance, (smooth) convolution
algebras of \'etale groupoids have local units, and this implies
$H$-unitality; actually we will show that
$C_{c}^{\infty}\bigl(\G\bigr)$ is $H$-unital for any smooth groupoid $\G$.

In the present context, the algebra $A$ we work with 
is endowed with
a locally convex topology, and the relevant homology groups are
obtained by replacing the algebraic tensor products by topological
ones. One has many topological tensor products available, but the
appropriate choice is often dictated by the type of algebras under
consideration and by the desire to have a computable target for Chern
characters. For instance, when $A=
C_{c}^{\infty}\bigl(M\bigr)$ for a manifold $M$, one recovers (compactly
supported) DeRham cohomology and the classical Chern character,
provided one uses the inductive tensor product of locally convex
algebras. The same product is relevant for convolution
algebras, and, in the sequel, $\otimes$ will denote this
topological tensor product. Actually, the only thing the reader
needs to know about it is that $C_{c}^{\infty}\bigl(M\bigr)\otimes
C_{c}^{\infty}\bigl(N\bigr)\cong C_{c}^{\infty}\bigl(M\times N\bigr)$ 
for any two
manifolds $M$, $N$ (and our results apply to any tensor product
with this property).

\section{Main results}
\label{sec2}

In this section we present our main results
concerning smooth groupoids which
appear in foliation theory. The first one is the characterisation
theorem already mentioned in the introduction:

\begin{st}\label{charact}For a smooth groupoid $G$, the following
are equivalent:
\begin{enumerate}[(i)]
\item
$G$ is Morita equivalent to a smooth \'etale groupoid;
\item
The Lie algebroid $\mathfrak{g}$ of $G$
has an injective anchor map;
\item
All isotropy Lie groups of $G$ are discrete.
\end{enumerate}
\end{st}

We will refer to groupoids  with this property as {\it foliation
groupoids}. For instance, the action groupoid $M\rtimes G$ associated 
to the action of
a Lie group on a manifold $M$ (which models the orbit space $M/G$) is 
a foliation
groupoid, provided all the isotropy groups $G_x= \{ g\in G: xg= x\}$ 
are discrete.
Also, if $\G$ is a foliation groupoid, then so is any pull-back of $\G$
(e.g. the groupoid $\G_{\U}$ associated to any cover $\U$ of $\nG{0}$).
The motivating examples are, however, the holonomy and the
monodromy  groupoids  $\Hol(M, \F)$ and $\Mon(M, \F)$ of any foliation 
$(M, \F)$ (note that the monodromy
groupoid appears in literature also under the name of ``the homotopy 
groupoid'' \cite{Mon}).
The construction of the holonomy along longitudinal paths (paths
inside leaves) can be viewed as a morphism
\begin{equation}\label{hol} \hol: \Mon(M, \F)\rmap \Hol(M, \F)
\end{equation}
which is the identity on $M$ (i.e., it is a morphism of groupoids over $M)$.

Note that any foliation groupoid $\G$ defines a foliation
$\F$ on $\nG{0}$, and $\G$ can be viewed as an integration of $\F$.
In many examples one actually starts with a foliation $(M, \F)$,
and then chooses a convenient foliation groupoid $\G$ integrating $\F$.
It is generally accepted that the holonomy and the monodromy groupoids are
actually extreme examples of such integrations.
The following proposition gives a precise formulation of this 
principle. For simplicity we restrict
ourselves to {\it $s$-connected} groupoids, i.e. groupoids $\G$ with 
the property that
all its $s$-fibers are connected. Recall \cite{McK} that, if $\G$ is 
arbitrary, one can allways find an open
$s$-connected subgroupoid of $\G$ by taking the connected components of the
units in the $s$-fibers.

\begin{prop}\label{classific}
Let $(M, \F)$ be a foliation. For any $s$-connected smooth groupoid 
$\G$ integrating
$\F$, there is a natural factorization of the holonomy morphism 
(\ref{hol}) into
homomorphisms $h_{\G}$, $\hol_{\G}$ of groupoids over $M$,
\[ \xymatrix{  \Mon(M, \F) \ar[r]^-{h_{\G}} & \G \ar[r]^-{\hol_{\G}} &
\Hol(M, \F) .    }\]
The maps $h_{\G}$ and $\hol_{\G}$ are surjective local 
diffeomorphisms. Moreover,
$\G$ is $s$-simply connected (i.e. has simply connected $s$-fibers) if and only
if $h_{\G}$ is an isomorphism.
\end{prop}

We will give explicit constructions of $h_{\G}$ and 
$\hol_{\G}$ later.
However, we should remark that the first of these homomorphisms is a
consequence of integrability results for Lie algebroids in \cite{McXu}; see
also \cite{MM}.

We next turn to the cyclic homology of convolution algebras of
foliation groupoids. Since the \'etale
case is well understood \cite{BrNi, Cra, CrMo}, our aim is to show that
the homology doesn't change when one passes from a given foliation
groupoid to a Morita equivalent \'etale groupoid. Thus, one of our main
results is the following:

\begin{st}\label{main2} If $\G$ and $\Ha$ are Morita equivalent
foliation groupoids, then
\[ HC_{*}\bigl(C_{c}^{\infty}(\G)\bigr)\cong
 HC_{*}\bigl(C_{c}^{\infty}(\Ha)\bigr)\ ,\]
and similarly for Hochschild and periodic cyclic homology.
\end{st}

We emphasize that, due to the applications we have
in mind, our aim is to prove the previous theorem by means of
explicit formulas (see the remark at the end). As said in the 
introduction, we conjecture that
this theorem in fact holds for smooth groupoids generally. Note
also that some of our lemmas are proved in this generality. For instance,
since $H$-unitality is usually relevant to excision theorems \cite{Wod},
and since convolution algebras appear in the short exact sequences
given by the pseudo-differential calculus \cite{WeNi},
the following result which is independent interest:

\begin{prop}\label{Hunital} The convolution algebra
$C_{c}^{\infty}\bigl(\G\bigr)$ of any smooth groupoid $\G$ is $H$-unital.
\end{prop}

Note that Theorem \ref{main2}, combined with the 
results of \cite{BrNi, Cra, CrMo}
concludes the computation of the cyclic homology for various 
foliation groupoids. Apart from the holonomy
and the monodromy groupoids, we mention the groupoids modeling 
orbifolds, and the groupoids
associated to Lie group actions with discrete stabilizers. \\

\begin{rmks}\rm
Before turning to the proofs in the next 
section, we make some further remarks:

\medskip
(i)
The holonomy groupoid of a foliation $(M, \F)$ appears as
the right model for the leaf space $M/\F$. Proposition \ref{classific} shows
that  it is the minimal {\it smooth} ``desingularization''
of the leaf space. We want to point out, however, that the holonomy 
groupoid may not be
the most appropriate model
when looking at problems whose primarly interest is not the leaf space.
One can find many examples where other foliation groupoids integrating $\F$
are equally good, and sometimes even more suitable. This applies, for example,
to the results of \cite{HeLa} which can be obtained using any 
Hausdorff groupoid
integrating the given foliation (all that matters is that the 
groupoid has the property stated in Lemma \ref{lemfolgr}
below). Regarding the Hausdorffness, we remark that
there is no relation between the Hausdorffness of $\Mon(M, \F)$ and of
$\Hol(M, \F)$, and there are foliations $\F$ whose monodromy and
holonomy groupoids are both non-Hausdorff, but which admit Hausdorff
integrations $\G$.

\medskip
(ii)
In the longitudinal index theory for foliations 
$(M, \F)$ of a compact manifold $M$, the
analytic index of a longitudinal elliptic operator $D$ can again be 
defined using any
foliation groupoid $\G$ integrating $\F$. First of all one lifts $D$ to
an operator along the $s$-fibers of $\G$, and then the 
pseudodifferential calculus
on $\G$ (namely the short exact sequence given by the symbol map of 
Theorem $8$ in \cite{WeNi},
and the boundary map of the long exact sequence it induces in 
$K$-theory) gives a precise
meaning to the index 
$\Ind_{\G}(D)\in K_{0}\bigl(C_{c}^{\infty}(\G)\bigr)$ 
depending just on
the symbol of $D$ (actually just on the induced class in 
$K^1\bigl(S^*\F\bigr)$). Classically,
this construction is applied to the holonomy groupoid, but Theorem 
\ref{classific}
shows that the best choice is the monodromy groupoid of $(M, \F)$, 
where $\Ind_{\G}(D)$
provides the maximal information. Since the monodromy groupoid of the 
foliation by one leaf
is (Morita equivalent to) the fundamental group of $M$, our remark 
agrees also with the
framework of the $L^2$-index theorem of Atiyah \cite{At} and the 
higher versions of Connes
and Moscovici \cite{CoMo} (see also \cite{Ni2}).

Now, the general Chern character in cyclic homology \cite{Co3},
combined with our Theorem \ref{main2}, and with the computations at 
units given in
Theorem 4.1.3. of \cite{Cra}, give a Chern character localized at units
$\Chern^1: K_{0}\bigl(C_{c}^{\infty}(\G)\bigr)\rmap H^{*}_{c}\bigl(\G\bigr)$ 
(in order to restrict
to units, we do have to assume $\G$ to be Hausdorff) .
The cohomology groups $H^{*}_{c}(\G)$ are the re-indexed
homology groups of \cite{CrMo} applied to any etale groupoid 
equivalent to $\G$.
The longitudinal index formula for foliations (non-commutative 
approach) gives a
topological interpretation for $\Chern^1\bigl(\Ind_{\G}(D)\bigr)$. More general
formulas should correspond to other localizations (cf 4.1.2 in \cite{Cra})
of the Chern character.

\medskip
(iii)
Following a different route (in the spirit of 
Bismut's approach to the
families index theorem), Heitsch--Lazarov \cite{HeLa} define certain 
cohomology classes
$\overline{\Chern}_{\mathcal{E}}(D) \in H^{*}_{c,\basic}(M/\F)$ playing the role of
``the Chern character of the index bundle''. Here $H_{c,\basic}^{*}(M/\F)$ are
the basic cohomology groups of Haefliger \cite{minimal}. The connection with
Connes approach (conjectured in \cite{HeLa}) can be described as follows.
For any integration $\G$ of $\F$ there is a tautological
map $j_{b}: H_{c}^{*}(\G)\rmap H_{c,\basic}^{*}(M/\F)$, which combined 
with $\Chern^1$ previously
described, induces a basic Chern character at units $\Chern_{\basic}^{1}:
K_{0}\bigl(C_{c}^{\infty}(\G)\bigr)\rmap H^{*}_{c,\basic}(M/\F)$. For a 
longitudinal elliptic
operator $D$ one gets 
$\Chern_{\basic}^{1}\bigl(\Ind_{\G}(D)\bigr)\in H^{*}_{c,\basic}(M/\F)$
independent of the choice of the Hausdorff integration $\G$.
Comparing the two longitudinal index theorems of \cite{Co3} and 
\cite{HeLa}, one sees that
(with the proper normalizations) 
$\Chern_{\basic}^{1}\bigl(\Ind_{\G}(D)\bigr)= 
\overline{\Chern}_{\mathcal{E}}(D)$.
Of course, an interesting question is to give a direct argument for 
this equality
between the basic Chern character of the analytical index, and the 
Chern character of the index bundle.
In this context we remark that, in contrast with $\Chern^1$, it is 
possible to describe
the basic Chern character $\Chern_{\basic}^{1}$ by relatively simple 
explicit formulas (with the help of connections),
using Haefliger's integration \cite{minimal} along leaves
and the non-commutative version \cite{Ka} of the Chern--Weil 
construction (see \cite{Crath} for details).
\end{rmks}

\section{Proof of the characterisation theorem}
\label{sec3}
%
In this section we present the proofs of Theorem \ref{charact}
and of Proposition \ref{classific}.

\paragraph{Proof of Theorem \ref{charact}:}~\\
$(ii)\Leftrightarrow (iii)$: 
immediate because the Lie algebra of the
isotropy group $G_{x}^{x}$ is the kernel of the anchor map
$\alpha:\mathfrak{g}_x \rightarrow T_x\bigl(G_0\bigr)$.

\noindent
$(i)\Rightarrow (iii)$: 
since the isotropy groups of
an \'etale groupoid are clearly discrete, it suffices to remark
that this property is invariant under Morita equivalence. Indeed,
since the pullback square (\ref{eq1}) has a surjective submersion
on the bottom, the fibers of the left hand vertical map are
discrete if and only if those of the right hand vertical map are.  
Thus, the
isotropy groups of $G$ are discrete precisely when those of $H$ are.

\noindent
$(ii) \Rightarrow (i)$: 
suppose the anchor map
\[  \alpha: \mathfrak{g} \rightarrow T\bigl(\nG{0}\bigr)\]
of the Lie algebroid $\mathfrak{g}$ of $G$ is injective.  Write
$\F\subseteq T\bigl(\nG{0}\bigr)$ for the
image $\alpha(\mathfrak{g})$.  Then $\F$ is an involutive subbundle
of $T\bigl(\nG{0}\bigr)$, hence defines
a foliation $\F$ of $\nG{0}$.  On the other hand, the submersion
$s: \nG{1}\rightarrow \nG{0}$ (source)
defines a foliation $\tilde{\F}$ on $\nG{1}$, whose leaves are the
connected components of the fibers of
$s: \nG{1} \rightarrow \nG{0}$. Denote by $p$ the dimension of $\F$,
by $q$ its codimension. From the hypothesis,
the dimension of $\tilde{\F}$ is $p= \dim(G_x)= \rank(\mathfrak{g})$,
while its codimension $n$ is equal to the codimension
of $G_x$ in $G$, so that $n= p+q$.

\begin{lem}\label{lcovproj} The target map $t: (\nG{1}, 
\tilde{\F})\rmap (\nG{0},
\F)$ maps leaves into leaves, and its restriction to each leaf is
a local diffeomorphism. If $\G$ is $s$-connected, then, for any
point $x\in \nG{0}$, the space $t(\G_x)= L_x$ is the leaf through
$x$, and
\begin{equation}\label{covproj} t: G_x \rmap L_x \end{equation}
is a smooth covering projection with structure group $G_{x}^{x}$.
\end{lem}

\paragraph{Proof:} For any $g: x\rmap y$ in $\G$, one has
  a commutative diagram
\[ \xymatrix{
T_{y}\bigl(G_y\bigr)= \mathfrak{g}_y
\ar[d]_-{(dR_g)_y}^-{\sim}\ar[rrd]^-{\alpha_y} &  
 & \\ T_{g}\bigl(G_x\bigr)=
\tilde{\F}_{g} \ar[rr]^-{(dt)_g} & & T_{y}\bigl(\nG{0}\bigr) }\] 
where
$\alpha_y$ maps $\mathfrak{g}_y$ isomorphically into $\F_y$, and
$R_{g}: (G_y, y)\rmap (G_{x}, g)$ is the right multiplication by
$g$. Thus, the target map induces an isomorphism
\begin{equation}\label{izt}
(dt)_g: \tilde{\F}_g
\stackrel{\sim}{\rmap} \F_y\;.
\end{equation}
This shows that the target map $t: (\nG{1}, \tilde{\F}) \rmap
(\nG{0}, \F)$ maps leaves to leaves,
  and that its restriction to each leaf is a local diffeomorphism. Hence,
for any $y\in \nG{0}$, and any connected component $C$ of
$\G_{y}$, the map $t |_{C}$ is a local diffeomorphism of $C$ into
some leaf of $\F$. To prove it is onto, it suffices to remark that
$\{ t(C): C$ is a connected component of $G_y, y\in \nG{0} \}$ is
a partition of $\nG{0}$. Indeed, if $C_i\subset G_{y_i}$ are
connected components so that $t(C_1)\cap t(C_2)$ is non-empty, we
find $g_i\in C_i$ with $t(g_1)= t(g_2)= y$. Since $R_{g_{i}}:
G_{y}\rmap G_{y_i}$ are diffeomorphisms, 
the $R_{g_{i}}^{-1}\bigl(C_i\bigr)$
will be connected components of $G_{y}$, both containing $1_y$,
hence $R_{g_{1}}^{-1}\bigl(C_1\bigr)= R_{g_{2}}^{-1}\bigl(C_2)$, 
which shows that $t(C_1)= t(C_2)$. 
\\\hglue .1cm\hfill\Boxe\par

The following Lemma will complete the proof of the Theorem:
\begin{lem}\label{mortrs} For any transversal $T$ of $\F$, the 
groupoid $G_{T}^{T}$
  is \'etale. If $T$ is complete, then $G$ is Morita equivalent to
$G_{T}^{T}$.
\end{lem}

\paragraph{Proof:} 
First we claim that the source map restricts to a
local diffeomorphism
\begin{equation}\label{locdifs} s: \G^{T}\rmap \nG{0} \ .
\end{equation}
Since $t$ is a submersion (hence, in particular, it is transversal
to $T$), $\G^T= t^{-1}(T)$ is a submanifold of $\nG{1}$ of
codimension equal to the codimension of $T$ in $\nG{0}$ (i.e. to
$p$), whose tangent space at $g: x\rmap y$ consists of vectors
$\xi\in T_{g}\bigl(\nG{1}\bigr)$ with the property that $(dt)_g(\xi)\in
T_{y}\bigl(T\bigr)$. By counting dimensions, it suffices to prove that the
map above is an immersion, i.e., since $\ker (ds)_{g}=
\tilde{\F}_{g}$, that
\begin{equation}\label{transint}
\tilde{\F}_{g}\cap (dt)_{g}^{-1}\bigl( T_{y}(T)\bigr)= \{ 0\} .
\end{equation}
But this is immediate from the isomorphism (\ref{izt}) and the
fact that $T$ is transversal to~$\F$.

Since
(\ref{locdifs}) is a local diffeomorphism, the inverse image
$G_{T}^{T}$ of $T$ is a submanifold, and the restriction $s:
\G_{T}^{T}\rmap T$ is a local diffeomorphism. Thus $\G_{T}^{T}$ is
\'etale. Moreover, if the transversal $T$ is complete, then $s:
\G^{T}\rmap \G$  is a surjection, and hence the obvious functor
$\G_{T}^{T}\rmap \G$ is an essential equivalence. This proves the lemma. 
\hfill\Boxe\par

For the proof of Proposition \ref{classific} we need 
the following Lemma. 
We first recall some terminology. Given a submersion $\pi: U\rmap T$,
the connected components of its fibers define a foliation
on $U$. Denote by $U\times_{T}U$ the fibered product
$\{ (x, y)\in U\times U: \pi(x)= \pi(y)\}$.
We say that $\pi$ is a  {\it trivializing submersion} of $\F$ if its
domain $U$ is open in $\nG{0}$, the fibers of $\pi$ are contractible, and
they coincide with the plaques of $\F$ in $U$.

\begin{lem}\label{lemfolgr} Let $\G$ be a foliation groupoid, and let
$\F$ be the induced foliation on $\nG{0}$. For any trivializing
submersion $\pi: U\rmap T$ of $\F$, there exists a unique open
subgroupoid $\G(U)\subset \nG{1}$ such that the map $(t, s):
\nG{1}\rmap \nG{0}\times \nG{0}$ restricts to an isomorphism of
smooth groupoids:
\[ (t, s): \G(U) \stackrel{\sim}{\rmap} U\times_{T}U \ .\]
\end{lem}

\paragraph{Proof:} First note that any such
open subgroupoid of $\G$ is contained in the $s$-connected component
of $G_{U}^{U}$. Hence it suffices to show that this $s$-connected
component, denoted $G(U)$, has the desired property. In other
words, it suffices to prove that if $U= \nG{0}$ and if $\G$ is
$s$-connected, the map $(t, s): \nG{1}\rmap U\times_{T}U$ is a
diffeomorphism. Remark that (\ref{izt}) implies that $(t, s):
\nG{1}\rmap U\times U$ is an immersion. By counting the
dimensions, it follows that $(t, s): \nG{1}\rmap U\times_{T} U$ is
a local diffeomorphism. It is also bijective because, by Lemma
\ref{lcovproj}, for any $x\in U$, the map $t: G_{x}\rmap
\pi^{-1}\bigl(\pi(x)\bigr)$ is a covering projection with connected total
space, and contractible base space, hence it is a diffeomorphism.
\\\hglue.1cm\hfill\Boxe\par

\paragraph{Proof of Proposition \ref{classific}:} 
Of course one can use
Lemma \ref{lcovproj} to define $h_{G}$. We indicate a slightly
different description, which immediately implies the smoothness of
$h_{G}$. Let $\alpha: [0, 1]\rmap L$ be a longitudinal path with
$\alpha(0)= x$, $\alpha(1)= y$. By the local triviality of $\F$
and the compactness of $\alpha[0, 1]$ we find a sequence $U_{i}$
of domains of trivializing submersions $\pi_{i}: U_i\rmap T_{i}$,
and real numbers $t_i$, so that:
\[ 0= t_0< t_1< \ .\ .\ .\ < t_{k}= 1,\ \ \alpha([ t_{i}, t_{i+1}])
\subset U_{i}, \ \ 0\leq i\leq k-1\ .\] From Lemma \ref{lemfolgr}
we find unique arrows $g_{i+1}: \alpha(t_{i}) \rmap
\alpha(t_{i+1})$ in $G(U_i)$; we put
\[ h_{G}(\alpha):\,= g_{k}\, g_{k-1}\ .\ .\ .\ g_{1} \in G \ .\]
This definition closely resembles the construction of the
holonomy, and, by the same arguments, $h_{G}(\alpha)$ depends just
on the homotopy class of $\alpha$. The smoothness of $h_{G}$ is
immediate now, since, near $\alpha$, the smooth structure of
$\Mon(M, \F)$ is defined precisely using such chains $\{ U_{i}\}$
covering $\alpha$. That $h_{\G}$ is surjective if $\G$ is
$s$-connected follows from the fact that on the $s$-fibers it is
precisely the projection $\tilde{L}_x\rmap G_x$ induced by the
covering projection of Lemma \ref{lcovproj}.

We now construct $\hol_{\G}$. Actually, since $h_{\G}$ is surjective
and we want $\hol_{\G}\compose h_{\G}= \hol$, we only have to show
that the holonomy class of $\alpha$ is determined by
$h_{\G}(\alpha)$. For this, we remark that the holonomy germ of
$\alpha$ can be defined directly in terms of the arrow
$h_{\G}(\alpha)$. More precisely, giving any arrow $g: x\rmap y$,
and any transversal $T$ containing $x$ and $y$, one obtains an
induced germ $\sigma_{g}^{T}: (T, x)\rmap (T, y)$ due to the fact
that $\G_{T}^{T}$ is \'etale (see our preliminaries on groupoids).
We claim that when $g= h_{\G}(\alpha)$, this germ $\sigma_{g}^{T}$
coincides with the holonomy germ of $\alpha$. This is clear when
$\alpha$ is contained in the domain of a trivializing submersion.
In general, we use that $\sigma_{g}^{T}$ is functorial in $g$, and
that $\sigma_{g}^{T}= \sigma_{g}^{S}$ whenever $S$ is another
transversal containing $T$. Choosing any transversal $T$
containing $x, y$ and all the $\alpha(t_i)\,$s above, it follows
that the germ associated to $g= h_{G}(\alpha)$, which is
\[ \sigma_{g}^{T}= \sigma_{g_k}^{T}\,
\sigma_{g_{k-1}}^{T}\ .\ .\ .\ \sigma_{g_1}^{T}: (T, x)\rmap (T,
y), \] coincides with the holonomy germ of $\alpha$. The last part
of the Corollary follows from Lemma \ref{lcovproj}. 
\hfill\Boxe\par

\section{Proof of the invariance theorem}
\label{sec4}
%
In this section we present the proofs of Theorem \ref{main2}
and of Proposition \ref{Hunital}. We will assume throughout that $\G$ 
is Hausdorff. However, we point out that
our proofs also apply to the non-Hausdorff case, provided one uses 
the compact supports defined in \cite{CrMoHa}
(similar extensions to the non-Hausdorff case already occur in 
\cite{Cra, CrMo}.)

\paragraph{Proof of Proposition \ref{Hunital}:} We first need some
remarks about the extension of compactly supported smooth
functions. Let $M$ be a manifold and let $A$ be a closed subset of
$M$. Write $\es_{M, A}$ for the fine sheaf of smooth functions on
$M$ which vanish on $A$. For a closed submanifold  $N\subset M$,
there is an obvious restriction
\begin{equation}\label{restr}
\Gamma_{c}(\es_{M, A})\rmap \Gamma_c(\es_{N, N\cap A})\ .
\end{equation}

We will be concerned with the surjectivity of this
map, for specific $M$, $N$, and $A$. Note that this surjectivity
is a local property: if each $x\in N$ has a neighborhood $U$ in
$M$ such that $\Gamma_{c}(\es_{U, A\cap U})$ $\rmap$
$\Gamma_c(\es_{N\cap U, N\cap A\cap U})$ is surjective, then it
follows (by a partition of unity or a Mayer--Vietoris argument)
that (\ref{restr}) is surjective.

If $x\in N$,
one can always choose a neighborhood $U$ of $x$ in $M$ and a
retraction $r: U\rmap U\cap N$. If, for any $x$, these $U$ and $r$ can
be chosen such that $r(A\cap U)\subset A$,  we say that $A$ is
{\it locally retractible} to $N$ in $M$. Note that this implies the
surjectivity of (\ref{restr}). Indeed, since both properties are
local, we may assume that there exists a retraction $r: M\rmap N$
such that $r(A)\subset A$; then, for any $\phi\in
C_{c}^{\infty}\bigl(N\bigr)$ vanishing on $A\cap N$, $\tilde{\phi}(x)=
\theta(x) \phi\bigl(r(x)\bigr)$ defines an extension of $\phi$ to $M$
vanishing on $A$, provided we choose 
$\theta\in C_{c}^{\infty}\bigl(M\bigr)$
with $\theta\equiv 1$ on the support of $\phi$.

An easy argument based on the canonical local form of a submersion
shows that:
\begin{lem}\label{retract} 
Let $s: X\rmap Z$ be a submersion, let $f: Y\rmap Z$
be a  smooth map, and let $B\subset Y$ be a closed subset. Then
$A= X\times B$ is locally retractible in $M= X\times Y$ to $N= X\times_{Z}
Y$. In particular, (\ref{restr}) is surjective.
\end{lem}

For the proof of the proposition, we have to prove
that if $\psi\in C_{c}^{\infty}\bigl(\nG{1}^{p}\bigr)$  is a cycle with
respect to $b'$ (i.e. $b'(\psi)= 0$), then it is
$b'$-homologous to zero (i.e. is of type $b'(\tilde{\psi})$
for some $\tilde{\psi}\in C_{c}^{\infty}\bigl(\nG{1}^{p+1}\bigr)$). 
Recall
that $b'(\psi)= \sum_{1}^{p-1} (-1)^{i} d_{i}(\psi)$, where
\[ d_i(\psi)(g_1, \ldots , g_{p-1})= \int_{uv= g_i} \psi(g_1,
\ldots , g_{i-1}, u, v, \ldots , g_{p-1}) \]
We will first show
that, for each $k= 1, \ldots p$, there exists a cycle $\psi_k$
homologous to $\psi$ such that
\begin{equation}\label{normalization}  \psi_k(g_1, \ldots , g_p)= 0\mbox{ if }
s(g_i)= t(g_{i+1}) \mbox{ for some } i< k .
\end{equation}
Notice that for such a cycle $\psi_k$, we have $d_{i}(\psi_k)= 0$
for $i< k$. We construct $\psi_k$ by induction on $k$. For $k=1$
the condition (\ref{normalization}) is vacuous, and we can take
$\psi_1= \psi$. Suppose $\psi_1, \ldots , \psi_k$ have been
defined. Let $K= \supp(\psi_k)$, which is a compact subset of
$\nG{1}^{p}$. 
Let $L= \{x: \exists (g_1, \ldots , g_p)\in K \mid x= t(g_k)\}$, 
and let 
$\theta\in C_{c}^{\infty}\bigl(\nG{1}\bigr)$ be
a function such that $\int_{t^{-1}(x)}\theta= 1$ for all $x\in L$.
Now define a function $\phi$ on the submanifold $N\subset M=
\nG{1}^{p+1}$ consisting of those $(g_1, \ldots , g_{p+1})$ for
which $g_kg_{k+1}$ is defined, by
\[ \phi(g_1, \ldots , g_{p+1})= \theta(g_k) \psi_k(g_1, \ldots ,
g_kg_{k+1}, \ldots , g_{p+1})\] Thus $\phi(g_1, \ldots , g_{p+1})=
0$ as soon as $s(g_i)= t(g_{i+1})$ for some $1\leq i < k$. By
Lemma \ref{retract}, we can find $\tilde{\phi}\in
C_{c}^{\infty}\bigl(\nG{1}^{p+1}\bigr)$ such that $\tilde{\phi}(g_1, \ldots
, g_{p+1})$ equals zero if $s(g_i)= t(g_{i+1})$ for some $1\leq i
< k$, and equals $\theta(g_k)\psi_k(g_1, \ldots , g_kg_{k+1},
\ldots , g_{p+1})$ if $s(g_k)= t(g_{k+1})$. Then we have
$d_i(\tilde{\phi})= 0$ for $1\leq i< k$, $d_k(\tilde{\phi})=
\psi_k$, while for $i> k$ and $s(g_k)= t(g_{k+1})$,
$d_i(\tilde{\phi})(g_1, \ldots , g_p)= \theta(g_k)
d_{i-1}(\psi_k)(g_1, \ldots , g_kg_{k+1}, \ldots , g_p)$. So,
still assuming $s(g_k)= t(g_{k+1})$,
\begin{eqnarray}
b'(\tilde{\phi})(g_1, \ldots , g_p)
  & = & (-1)^{k}\psi_k(g_1,
\ldots , g_p) + \sum_{j=k}^{p} (-1)^{j+1}\theta(g_k)
d_{j}(\psi_k)(g_1, \ldots , g_kg_{k+1}, \ldots , g_p) \nonumber\\
  & = & (-1)^k \psi_k(g_1, \ldots , g_p) - \theta(g_k)
b'(\psi_k)(g_1, \ldots , g_kg_{k+1}, \ldots , g_p) \nonumber \\
  & = & (-1)^k \psi_k(g_1, \ldots , g_p)
\end{eqnarray}
Thus we can put $\psi_{k+1}= \psi_k- (-1)^k b'(\tilde{\phi})$ to
obtain the desired property. Having thus defined $\psi_1, \ldots ,
\psi_p$, the construction of $\tilde{\phi}$ for $k= p$ gives a
function $\tilde{\phi}$ with $d_i(\tilde{\phi})= 0$ for $i< p$,
and $d_p(\tilde{\phi})= \psi_p$. Thus $b'(\tilde{\phi})=
(-1)^{p} \psi_p$, showing that $\psi_p$ is a boundary. This proves 
Proposition \ref{Hunital}.\hfill $\Boxe$\par

To prove Theorem \ref{main2} we need some preliminary lemmas.
We first compare the convolution algebra of $\G$, with
the one of the groupoid $\G_{\U}$ induced by $\G$ and an open
covering $\U= \{ U_{i}\}$ of $\nG{0}$ (cf. our preliminaries). The
elements of $ C_{c}^{\infty}\bigl(\G_{\U}\bigr)= \oplus_{i, j}
C_{c}^{\infty}\bigl(\G_{U_j}^{U_i}\bigr)$  can be written as matrices $\f=
(\f_{i, j})_{i, j}$. We will also use the following left/right
action of $C_{c}^{\infty}\bigl(\nG{0}\bigr)$ 
on $C_{c}^{\infty}\bigl(\G\bigr)$:
\[ (f \phi) (g)= f(t\ps g\pd) \phi ( g ),\ \ (\phi f) (g)= \phi ( g ) f
(s\ps g \pd)\ .\]


\begin{lem}\label{redcov} For any smooth groupoid $\G$, and for
any family $\{ \lambda_i\in C_{c}^{\infty}\bigl(\nG{0}\bigr)\}$ so that
$\{\lambda_{i}^{2}\}$ is a partition of unity subordinated to a
locally finite open covering $\U= \{ U_i\}$ of $\nG{0}$, the map
\begin{equation}
\label{lambda} \lambda: C_{c}^{\infty}\bigl(\G\bigr) \rmap
C_{c}^{\infty}\bigl(\G_{\U}\bigr), \ \lambda(\phi)=
(\lambda_{i}\phi\lambda_{j})_{i, j} .
\end{equation}
is an algebra homomorphism whose induced maps in cyclic type
homologies are injective.
\end{lem}

\paragraph{Proof:} One has an obvious inclusion $i:
C_{c}^{\infty}\bigl(\G_{\U}\bigr)\hookrightarrow
M_{\infty}\bigl(C_{c}^{\infty}(\G)\bigr)$, already suggested by the notation
for elements of $C_{c}^{\infty}\bigl(\G_{\U}\bigr)$.  It suffices to prove
that the composition $i^{\, \lambda}= i\compose \lambda$ induces
isomorphism in cyclic homologies. Using the algebras $A=
C_{c}^{\infty}\bigl(\G\bigr)$, $A_0= C_{c}^{\infty}\bigl(\nG{0}\bigr)$, 
we find
ourselves in the abstract situation where we have a triple
\[ (A, A_0, \lambda), \ \  \ \lambda= \{\lambda_i\} \]
where $A$ is an $H$-unital algebra, $A_0$ is an algebra which acts
on both sides on $A$, and the $\lambda_i\in A_0$ are elements so
that, for any given $a\in A$ or $a\in A_0$, the products
$\lambda_ia$ and $a\lambda_i$ are nonzero for only finitely many
$i$, and $\sum\lambda_{i}^{2}a$ $=$ $\sum a\lambda_{i}^2= a$. We
prove that, in this situation, the algebra homomorphism $i^{\,
\lambda}: A\rmap M_{\infty}(A)$, $i^{\, \lambda}(a)=
(\lambda_{i}a\lambda_{j})_{i, j}$ induces isomorphisms in the
cyclic type homologies.

Let us first consider the special case where $A_0$ is a 
subalgebra of $A$. Recall from \cite{Lo} that the trace map $\Trace_{*}:
C_*\bigl(M_{\infty}(A)\bigr)\rmap C_*(A)$,
\[ \Trace_*(a^0, a^1, \, . \, .\, .\, , a^n)= \sum_{i_0, . . . , i_n}
(a^{0}_{i_0i_1}, a^{1}_{i_1i_2}, \, . \, .\, .\, ,
a^{n}_{i_ni_0}),\ a^{i}\in M_{\infty}(A),\] has this property.
Using the $SBI$ argument, it now suffices to show that $\Trace_{*}
i^{\, \lambda}_{*}$, acting on the Hochschild complex 
$\bigl(C_*(A), b\bigr)$, 
is homotopic to the
identity. For this, we construct the homotopy
\begin{eqnarray}
h(a^0, a^1, \ldots , a^n) & = & \sum_{i_0, . . . , i_n}
(a^0\lambda_{i_0}, \lambda_{i_0}a^1\lambda_{i_1}, \ldots ,
\lambda_{i_{n-1}}a^{n}\lambda_{i_n}, \lambda_{i_n})-   \nonumber\\
  &- & \sum_{i_0, . . . , i_{n-1}} (a^0\lambda_{i_0},
\lambda_{i_0}a^1\lambda_{i_1}, \ldots ,
\lambda_{i_{n-2}}a^{n-1}\lambda_{i_{n-1}}, \lambda_{i_{n-1}}, a^n)
+ \ldots +                                             \nonumber\\
  & + & (-1)^{n-1} \sum_{i_0, i_1} (a^0\lambda_{i_0}, 
\lambda_{i_0}a^1\lambda_{i_1},
\lambda_{i_1}, a^2, \ldots , a^n)+                     \nonumber\\
  & + & (-1)^n\sum_{i_0} (a^0\lambda_{i_0}, \lambda_{i_0}, a^1, \ldots , a^n)
.\nonumber
\end{eqnarray}

In the general case, we use the new algebra
$A_0\cross A$ which is $A_0\oplus A$ with the product
\[ (\lambda, a)(\eta, b)= (\lambda\eta, \lambda b+ a\eta+ ab), \
\lambda, \eta\in A_0, a, b\in A ,\] Remark that $A_0\cross A$
contains $A_0$ as a subalgebra (with the inclusion $\rho(\lambda)=
(\lambda, 0)$), and the map $i^{\, \lambda}$ lifts to a map
between short-exact sequences
\[ \xymatrix{
0\ar[r] & A\ar[r]^-{i}\ar[d]^-{i^{\, \lambda}} &  A_0\cross
A\ar[r]^-{\pi}\ar[d]^-{i^{\, \lambda}}&  A_0\ar[d]^-{i^{\,
\lambda}}\ar[r] & 0 \\ 0\ar[r] & M_{\infty}(A)\ar[r]^-{\tilde{i}}
&  M_{\infty}(A_0\cross A)\ar[r]^-{\tilde{\pi}}&  M_{\infty}(A_0)
\ar[r] & 0 }\] ($i(a)= (0, a)$, $\pi(\lambda, a)= \lambda$, and
$\tilde{i}, \tilde{\pi}$ are induced by $i$, $\pi$). Note that $\rho$ 
is an algebra splitting of $\pi$.
By the previous discussion, the statement is true for $A_0\cross
A$ and $A_0$; to deduce it for $A$, it suffices to use Wodzicki's
excision \cite{Wod} for Hochschild/cyclic homology. \hfill$\Boxe$\par

Next, let $R$ be the pair groupoid
$\mathbb{R}^p\times \mathbb{R}^p$ over $\mathbb{R}^p$. For
foliation groupoids $\G$ we will see that, for suitable choices of
coverings $\U$, the groupoid $\G_{\U}$ becomes isomorphic to $\Gamma\times R$
for some \'etale groupoid $\Gamma$. Therefore, we state and prove the 
following lemma only for such groupoids
$\Gamma$ (we mention however that, using Lemma
\ref{Hunital}, one can actually prove it for general smooth
groupoids).

\begin{lem}\label{prodet} For any \'etale groupoid $\Gamma$, there is 
an isomorphism
\[ \tau_{*}: HC_{*}\bigl( C_{c}^{\infty}(\Gamma\times R)\bigr) 
\stackrel{\sim}{\rmap}
HC_{*}\bigl( C_{c}^{\infty}(\Gamma)\bigr), \]
and similarly for Hochschild
and periodic cyclic homology (see below for explicit formulas).
\end{lem}

\paragraph{Proof:} 
The convolution algebra $\R$ of $R$ consists of
compactly supported smooth functions on $k(x, y)$ on
$\mathbb{R}^p\times \mathbb{R}^p$, with the product
\[ (k_1k_2)(x, z)= \int k_1(x, y)k_2(y, z) dy\ .\]
One has the usual trace $\tau$ on $\R$,
\begin{equation}\label{thetrace}
\tau: \R\rmap  \mathbb{C}, \ \tau(k)= \int k(x, x)dx ,
\end{equation}
and an induced chain map
\begin{eqnarray}
\tau_*: C_*\bigl(C_{c}^{\infty}(\Gamma\times R)\bigr)\rmap
C_*\bigl(C_{c}^{\infty}(\Gamma)\bigr), \nonumber \\
  \tau_*(a^0\otimes k^0, \, .\, .\, .\, , a^n\otimes k^n):= \tau(k^0\, 
.\, .\, .\, k^n)
  (a^0, \, .\, .\, .\, , a^n), \label{fortau}
\end{eqnarray}
We choose $u$, and then define $\alpha$ as in
\begin{equation}\label{chooseu}
u\in C_{c}^{\infty}\bigl(I^{p}\bigr), \ \ \int u(x)^2 dx= 1\ , \ \alpha : =
u\otimes u\in \R
\end{equation}
(where $I^{p}= (-1, 1)^p$) and consider the algebra homomorphism
\[ j_{\alpha}: C_{c}^{\infty}\bigl(\Gamma\bigr) \rmap 
C_{c}^{\infty}\bigl(\Gamma\times R\bigr), \ j_{\alpha}(a)= a\otimes \alpha \]
Since $\tau_{*}j_{\alpha}= \Id$, it suffices to show that
\[ j_{\alpha}\tau_{*}: C_*\bigl(C_{c}^{\infty}(\Gamma\times R)\bigr)\rmap
C_*\bigl(C_{c}^{\infty}(\Gamma\times R)\bigr), \]
\[ (a^0\otimes k^0, \, .\, .\, .\, , a^n\otimes k^n)\mapsto
\tau(k^0\, .\, .\, .\, k^n) (a^0\otimes\alpha , \, .\, .\, .\,
a^n\otimes\alpha) \]
induces the identity in Hochschild homology (hence,
by the usual $SBI$-argument, in all cyclic homologies). Let us
first assume that $\Gamma_{0}$ is compact. We then have the
following homotopy:
\begin{eqnarray*} h (a^0\otimes x^0\otimes y^0, a^1\otimes x^1\otimes 
y^1, \ .\ .\ .\ ,
a^n\otimes x^n\otimes y^n)  = \\
  =\sum_{k=0}^{n} (-1)^k \tau(y^0\otimes x^1) \tau(y^1\otimes x^2)  \, 
.\, .\, .\,
  \tau(y^{k-1}\otimes x^k) \psi_{k}\end{eqnarray*}
where $\psi_k$ is the element
\[ (a^0\otimes x^0\otimes u, a^1\otimes u\otimes u, \ .\ .\ .\ , 
a^k\otimes u\otimes u,
1\otimes u\otimes y^k, a^{k+1}\otimes x^{k+1}\otimes y^{k+1}, \ .\ .\ 
.\ , a^{n}\otimes
x^{n}\otimes y^{n}) \]
for all $a^i\in C_{c}^{\infty}\bigl(\Gamma\bigr)$, 
$x^i\otimes y^i\in \R$. It is straightforward
to write the corresponding formula for the general elements in
$C_{c}^{\infty}\bigl(\Gamma\times \mathbb{R}^p\times \mathbb{R}^p\bigr)$).
When $\Gamma_{0}$ is not compact, we have to replace the unit
$1\in C_{c}^{\infty}\bigl(\Gamma_{0}\bigr)\subset 
C_{c}^{\infty}\bigl(\Gamma\bigr)$
appearing in the previous formula, by local units (compactly supported
smooth functions on $\Gamma_{0}$, which are constantly $1$ on
compacts which exhaust $\Gamma_{0}$.\hfill\Boxe\par

\paragraph{Proof of Theorem \ref{main2}:} 
Since the theorem 
is known for \'etale groupoids
\cite{Cra, CrMo}, and since any foliation groupoid $\G$ is Morita 
equivalent to an etale one
(e.g. $\G_{T}^{T}$ of Lemma \ref{mortrs}),
it suffices to find, for a given foliation groupoid $\G$, a complete 
transversal $T$ for which we can prove
that $HC_{*}\bigl(C_{c}^{\infty}(\G)\bigr)\cong 
HC_{*}\bigl(C_{c}^{\infty}(\G_{T}^{T})\bigr)$.
  Let $\U= \{ U_1, U_2, \ldots \}$ be a locally
finite cover of $\nG{0}$ by foliation charts, say $\f_i:
\mathbb{R}^p\times \mathbb{R}^p\cong U_i$ $\subset \nG{0}$, and
write $T_i= \f_{i}(\{ 0\}\times \mathbb{R}^q)\subset U_i$ for the
transversals, and $\pi_i: U_i\rmap T_i$ for the evident
projections. Furthermore, let $\G_{\U}$ be the groupoid induced by
the cover $\U$ as described in the preliminaries. Now observe
that, by Lemma \ref{lemfolgr}, there are isomorphisms $H_i:
U_{i}\times_{T_i}U_{i} \stackrel{\sim}{\rmap} G(U_i)$. Each such
$H_i$, $H_j$ induce a map
\[ h_{i, j}: \G_{U_j}^{U_i} \rmap   \G_{T_j}^{T_i} \]
\begin{equation}\label{forhij} h_{i, j}(g)= H_{i}(\pi_i(t\ps g\pd ) , 
t\ps g\pd)\ \compose\
g\ \compose\  H_{j}( s\ps g\pd, \pi_{j}( s\ps g\pd )
.\end{equation} If we write $T= \coprod T_i$ for the complete
transversal and $R$ for the pair groupoid of Lemma \ref{prodet},
we then obtain an isomorphism (compare to \cite{HiSk})
\begin{equation}\label{has}
h: \G_{\U} \stackrel{\sim}{\rmap}   \G_{T}^{T}\times R ,\ \ \
h(i, g, j)= (h_{i, j}(g), p_i(t\ps g\pd ), p_j(s\ps g\pd ).
\end{equation}
which can be described in terms of the $h_{i, j}$ and the
projections $p_i: U_i\rmap \mathbb{R}^{p}$ on the first
coordinate, by \[ h(i, g, j)= (h_{i, j}(g), p_i(t\ps g\pd ),
p_j(s\ps g\pd ).\]
The isomorphism $h$, combined with the map
$j_{\alpha}$ of the proof of Lemma \ref{prodet}, gives a map
\begin{equation}\label{hjalpha}
hj_{\alpha}: C_{c}^{\infty}\bigl(\G_{T}^{T}\bigr)\rmap
C_{c}^{\infty}\bigl(\G_{\U}\bigr)
\end{equation} which induces isomorphisms in cyclic type
homologies. Now consider a
sequence of smooth functions $\lambda_i\in C_{c}^{\infty}\bigl(\G\bigr)$
such that the $\lambda_{i}^{2}$ form
a partition of unity subordinate to $\U$. We can choose the $U_i$
and $\lambda_i$ in such a way (see the proof of the preliminary lemma 
of \cite{HiSk})
that for the open sets
$V_i= \f_i(I^{p}\times I^{q})$ with transversals $S_i=
\f_i(\{0\}\times I^{q})$ (recall that $I= (-1, 1)$) one has that
$\bar{V}_{i}\cap \bar{V}_{j}= \emptyset$ whenever $i\neq j$, while
$\lambda_{i} |_{V_{i}}= 1$ and each leaf of $\F$ meets at least
one $S_i$.

There is an obvious analogue of
(\ref{hjalpha}) associated to the family $\nV$ and to the complete
transversal $S= \coprod S_i$, and we obtain a commutative square:
\[ \xymatrix{
C_{c}^{\infty}\bigl(\G_{S}^{S}\bigr) \ar[r]^-{hj_{\alpha}} \ar[d]_-{e'} &
C_{c}^{\infty}\bigl(\G_{\nV}\bigr) \ar[d]^-{e}\\
C_{c}^{\infty}\bigl(\G_{T}^{T}\bigr)\ar[r]^-{hj_{\alpha}} &
C_{c}^{\infty}\bigl(\G_{\U}\bigr) } \] where the vertical $e$ and $e'$ are
given by extension by zero. In this diagram, the maps $hj_{\alpha}$
have been shown to induce isomorphisms in cyclic type homologies,
while the map $e'$ does so by Morita invariance for etale
groupoids \cite{Cra, CrMo}. Hence the map $e$ also induces such
isomorphisms. We also have a commutative diagram:
\[ \xymatrix{
C_{c}^{\infty}\bigl(\G_{\U}\bigr) \ar[r]^-{e''} \ar[dr]_-{e} &
C_{c}^{\infty}\bigl(\G\bigr) \ar[d]^-{\lambda} \\
  & C_{c}^{\infty}\bigl(\G_{\U}\bigr) } \]
where $\lambda$ is the map defined in Lemma \ref{redcov}, and
$e''$ is again defined by extension by zero. This diagram and
the previous remark on $e$ imply that the maps induced by
$\lambda$ in the cyclic homologies are surjective. Using
Lemma \ref{redcov}, it then follows that all the maps in the last
two diagrams induce isomorphisms in the cyclic type homologies. 
\hfill $\Boxe$\par

\begin{rmk}\rm
Let $\G$ be a foliation groupoid, and let $T$, $S$ be
the complete transversals previously constructed. There is a
commutative diagram
\[ \def\objectstyle{\textstyle}
\xymatrix{
  HC_{*}\bigl(C_{c}^{\infty}(\G)\bigr) \ar[r]^-{B} &
  HC_{*}\bigl(C_{c}^{\infty}(\G_{T}^{T})\bigr) \\
   & HC_{*}\bigl(C_{c}^{\infty}(\G_{S}^{S})\bigr)
     \ar[u]_-{e}\ar[lu]^-{A} } \]
where $A$, $B$, $e$ are isomorphisms described as follows:
\begin{enumerate}[(i)]
\item
$e$ is induced by the extension by zero map;
\item
$A$ is induced by the algebra homomorphism $A:
C_{c}^{\infty}\bigl(G_{S}^{S}\bigr)\rmap C_{c}^{\infty}\bigl(\G\bigr)$,
\[
A(\phi_{i,j})(g) = u_{i}\bigl(t(g)\bigr) 
   \phi\bigl(h_{i, j}(g)\bigr) u_{j}\bigl(s(g)\bigr)
\]
for all $\phi_{i, j}\in C_{c}^{\infty}\bigl(\G_{S_i}^{S_j}\bigr)$. 
Here $h_{i,
j}$ is given by the formula (\ref{forhij}), and $u_i= u\compose
p_i\in C_{c}^{\infty}\bigl(\nG{0}\bigr)$, with $u$ chosen as in
(\ref{chooseu}), and $p_i: U_i\rmap \mathbb{R}^p$ the
projection;
\item
$B$ is induced by the composition
\[ C_{*}\bigl(C_{c}^{\infty}(\G)\bigr) \stackrel{\lambda}{\rmap}
C_{*}\bigl(C_{c}^{\infty}(\G_{\U})\bigr) \stackrel{h}{\rmap}
C_{*}\bigl(C_{c}^{\infty}(\G_{T}^{T}\times R)\bigr) \stackrel{\tau_*}{\rmap}
C_{*}\bigl(C_{c}^{\infty}(G_{T}^{T})\bigr)\] where $\lambda$ is given by
(\ref{lambda}), $h$ is induced by the isomorphism (\ref{has}), and
$\tau_{*}$ is given by the formula (\ref{fortau}).
\end{enumerate}
\end{rmk}

\end{document}